% begin.tex

%=========================================================================%
% load begin.tex only once, but keep count to match \bye commands
%=========================================================================%

\ifx\begin\undefined\else\global\advance\srcdepth by
1\expandafter \fi

\def\begin{}
\newcount\srcdepth
\srcdepth=1

\outer\def\bye{\global\advance\srcdepth by -1
  \ifnum\srcdepth=0
    \def\endcmd{\vfill\eject\nopagenumbers\par\vfill\supereject\end}
  \else\def\endcmd{}\fi
  \endcmd
}

%=========================================================================%
% initialize TeX
%=========================================================================%

%\magnification=\magstephalf
%\magnification=\magstep 1
\baselineskip=13pt
%baselineskip=12pt
\hsize = 5.5truein
\hoffset = 0.5truein
\vsize = 8.5truein
\voffset = 0.2truein
\emergencystretch = 0.05\hsize

\overfullrule=0pt

\newif\ifblackboardbold

% comment out the following line if AMS msbm fonts aren't available
\blackboardboldtrue

%=========================================================================%
% select fonts
%=========================================================================%

\font\sectionfont=cmbx12

% Establish AMS blackboard bold fonts without using amssym.def, amssym.tex

\newfam\bboldfam
\ifblackboardbold
\font\tenbbold=msbm10
\font\sevenbbold=msbm7
\font\fivebbold=msbm5
\textfont\bboldfam=\tenbbold
\scriptfont\bboldfam=\sevenbbold
\scriptscriptfont\bboldfam=\fivebbold
\def\bbold{\fam\bboldfam\tenbbold}
\else
\def\bbold{\bf}
\fi

%=========================================================================%
% font size-changing command ("A Beginner's Book of TeX" p35, p275)
%=========================================================================%

\font\Arm=cmr8
\font\Ai=cmmi8
\font\Asy=cmsy8
\font\Abf=cmbx8
\font\Brm=cmr6
\font\Bi=cmmi6
\font\Bsy=cmsy6
\font\Bbf=cmbx6
\font\Crm=cmr5
\font\Ci=cmmi5
\font\Csy=cmsy5
\font\Cbf=cmbx5

\ifblackboardbold
\font\Abbold=msbm10 at 8pt
\font\Bbbold=msbm7 at 6pt
\font\Cbbold=msbm5
\fi

\def\smallmath{%
\textfont0=\Arm \scriptfont0=\Brm \scriptscriptfont0=\Crm
\textfont1=\Ai \scriptfont1=\Bi \scriptscriptfont1=\Ci
\textfont2=\Asy \scriptfont2=\Bsy \scriptscriptfont2=\Csy
\textfont\bffam=\Abf \scriptfont\bffam=\Bbf \scriptscriptfont\bffam=\Cbf
\def\rm{\fam0\Arm}\def\mit{\fam1}\def\oldstyle{\fam1\Ai}%
\def\bf{\fam\bffam\Abf}%
\ifblackboardbold
\textfont\bboldfam=\Abbold
\scriptfont\bboldfam=\Bbbold
\scriptscriptfont\bboldfam=\Cbbold
\def\bbold{\fam\bboldfam\Abbold}%
\fi
}

%=========================================================================%
% single-pass symbolic theorem labeling
%=========================================================================%

% Because this is a single-pass mechanism with no .aux file, forward
% references need to be declared in advance:

%   \forward{thm:main}{Theorem}{1.1}

% This is also the mechanism for "timely" declaration of labels, which
% will usually be buried within the corresponding theorem macros.
% A warning is issued if a label redeclaration is inconsistent, allowing
% forward references to be manually fixed.

%   \ref{thm:main} produces "Theorem~1.1"
%   \refs{thm:main} produces "Theorems~1.1"
%   \refn{thm:main} produces "1.1"

% Some TeX adapted from "The Advanced TeXbook" by David Salomon, chapter 9.

% Implementers: The code for \forward is subtle. Its second argument must
% be provided literally, e.g. "Theorem" rather that "\capitalize{theorem}".
% Its third argument must either be literal or a macro that expands
% directly to a literal, e.g. "\edef\numtoks{\number\proccount}".
% This use of \edef cannot be replaced by \def, which defers expansion.
% Failure to follow these rules will cause spurious warnings that forward
% references are inconsistent, when they are in fact consistent after
% expansion. Note the "Towers of Palo Alto" recreational math problem
% involving the iterated use of \expandafter to expand the first argument
% to \forwardsub before calling it.

\newlinechar=`@
\def\forwardmsg#1#2#3{\immediate\write16{@*!*!*!* forward reference should
be: @\noexpand\forward{#1}{#2}{#3}@}}
\def\nodefmsg#1{\immediate\write16{@*!*!*!* #1 is an undefined reference@}}

\def\forwardsub#1#2{\def\newref{{#2}{#1}}}

\def\forward#1#2#3{%
\expandafter\expandafter\expandafter\forwardsub\expandafter{#3}{#2}
\expandafter\ifx\csname#1\endcsname\relax\else%
\expandafter\ifx\csname#1\endcsname\newref\else%
\forwardmsg{#1}{#2}{#3}\fi\fi%
\expandafter\let\csname#1\endcsname\newref}

\def\firstarg#1{\expandafter\argone #1}\def\argone#1#2{#1}
\def\secondarg#1{\expandafter\argtwo #1}\def\argtwo#1#2{#2}

\def\ref#1{\expandafter\ifx\csname#1\endcsname\relax
  {\nodefmsg{#1}\bf`#1'}\else
  \expandafter\firstarg\csname#1\endcsname
  ~\expandafter\secondarg\csname#1\endcsname\fi}

\def\refs#1{\expandafter\ifx\csname#1\endcsname\relax
  {\nodefmsg{#1}\bf`#1'}\else
  \expandafter\firstarg\csname #1\endcsname
  s~\expandafter\secondarg\csname#1\endcsname\fi}

\def\refn#1{\expandafter\ifx\csname#1\endcsname\relax
  {\nodefmsg{#1}\bf`#1'}\else
  \expandafter\secondarg\csname #1\endcsname\fi}

%=========================================================================%
% widow control
%=========================================================================%

% usage:
% \widow{.2} % start new page if <.2 page left

\def\widow#1{\vskip 0pt plus#1\vsize\goodbreak\vskip 0pt plus-#1\vsize}

%=========================================================================%
% sections and theorems
%=========================================================================%

% use \showlabels or \showlabelsabove to display section and theorem labels

\def\marginlabel#1{}

\def\showlabelsabove{
\font\labelfont=cmss10 at 6pt
\def\marginlabel##1{\rlap{\smash{\raise 10pt\hbox{\labelfont##1}}}}
}

\newcount\seccount
\newcount\proccount
\seccount=0
\proccount=0

\def\stdskip{\vskip 9pt plus3pt minus 3pt}
\def\stdbreak{\par\removelastskip\penalty-100\stdskip}

\def\proof{\stdbreak\noindent{\sl Proof. }}

\def\qed{\vrule height 1.2ex width .9ex depth .1ex}

\def\Box{
  \ifmmode\eqno\qed
  \else\ifvmode\removelastskip\line{\hfil\qed}
  \else\unskip\quad\hskip-\hsize
    \hbox{}\hskip\hsize minus 1em\qed\par
  \fi\stdbreak\fi}

\def\references{
  \removelastskip
  \widow{.05}
  \vskip 24pt plus 6pt minus 6 pt
  \leftline{\sectionfont References}
  \nobreak\stdskip\noindent}

\def\ifempty#1#2\endB{\ifx#1\endA}
\def\makeref#1#2#3{\ifempty#1\endA\endB\else\forward{#1}{#2}{#3}\fi}

\outer\def\section#1 #2\par{
  \removelastskip
  \global\advance\seccount by 1
  \global\proccount=0\relax
                \edef\numtoks{\number\seccount}
  \makeref{#1}{Section}{\numtoks}
  \widow{.05}
  \vskip 24pt plus 6pt minus 6 pt
  \message{#2}
  \leftline{\marginlabel{#1}\sectionfont\numtoks\quad #2}
  \nobreak\stdskip}

\def\proclamation#1#2{
  \outer\expandafter\def\csname#1\endcsname##1 ##2\par{
  \stdbreak
  \advance\proccount by 1
  \edef\numtoks{\number\seccount.\number\proccount}
  \makeref{##1}{#2}{\numtoks}
  \noindent{\marginlabel{##1}\bf #2 \numtoks\enspace}
  {\sl##2\par}
  \stdbreak}}

\def\othernumbered#1#2{
  \outer\expandafter\def\csname#1\endcsname##1{
  \stdbreak
  \advance\proccount by 1
  \edef\numtoks{\number\seccount.\number\proccount}
  \makeref{##1}{#2}{\numtoks}
  \noindent{\marginlabel{##1}\bf #2 \numtoks\enspace}}}

\proclamation{definition}{Definition}
\proclamation{lemma}{Lemma}
\proclamation{proposition}{Proposition}
\proclamation{theorem}{Theorem}
\proclamation{corollary}{Corollary}
\proclamation{conjecture}{Conjecture}

\othernumbered{example}{Example}
\othernumbered{remark}{Remark}
\othernumbered{construction}{Construction}

%=========================================================================%
% enable postscript illustrations using epsf.tex
%=========================================================================%

% Usage:
% \draw{70}{fig}{} % draw fig.eps at 70% scale
% \draw{999}{fig}{} % draw fig.eps scaled to width of page

% Optional third argument can be multiple calls to \figtext; see below.
% More generally, the third argument is read in vertical mode, with the
% reference point at the lower left corner of the eps picture, whose
% dimensions are contained in the dimen registers \drawx and \drawy.
% This enables using TeX to generate the text that goes with the picture.
% To request that the picture be widened to respect the added text, 
% examine and modify the dimen registers \ngap, \egap, \sgap, \wgap.
% This is done automatically by the \figtext macro.

% These macros rely on "epsf.tex" which is the lowest level interface
% available for including encapsulated Postscript files in TeX documents.
% Rather that manually reading the .eps file to compute the nominal size,
% the \epsfbox macro is called twice, and two of its internal registers
% are examined after the first call. A major change to epsf.tex (unlikely)
% will require changes here. 

%\input epsf

\newcount\figcount
\figcount=0
\newbox\drawing
\newcount\drawbp
\newdimen\drawx
\newdimen\drawy
\newdimen\ngap
\newdimen\sgap
\newdimen\wgap
\newdimen\egap

\def\drawbox#1#2#3{\vbox{
  \setbox\drawing=\vbox{\offinterlineskip\epsfbox{#2.eps}\kern 0pt}
  \drawbp=\epsfurx
  \advance\drawbp by-\epsfllx\relax
  \multiply\drawbp by #1
  \divide\drawbp by 100
  \drawx=\drawbp truebp
  \ifdim\drawx>\hsize\drawx=\hsize\fi
  \epsfxsize=\drawx
  \setbox\drawing=\vbox{\offinterlineskip\epsfbox{#2.eps}\kern 0pt}
  \drawx=\wd\drawing
  \drawy=\ht\drawing
  \ngap=0pt \sgap=0pt \wgap=0pt \egap=0pt 
  \setbox0=\vbox{\offinterlineskip
    \box\drawing \ifgridlines\drawgrid\drawx\drawy\fi #3}
  \kern\ngap\hbox{\kern\wgap\box0\kern\egap}\kern\sgap}}

\def\draw#1#2#3{
  \setbox\drawing=\drawbox{#1}{#2}{#3}
  \advance\figcount by 1
  \goodbreak
  \midinsert
  \centerline{\ifgridlines\boxgrid\drawing\fi\box\drawing}
  \smallskip
  \vbox{\offinterlineskip
    \centerline{Figure~\number\figcount}
    \smash{\marginlabel{#2}}}
  \endinsert}

\def\nextfigtoks{%
  \advance\figcount by 1%
  \edef\numtoks{\number\figcount}%
  \advance\figcount by -1}

\newif\ifgridlines
\newbox\figtbox
\newbox\figgbox
\newdimen\figtx
\newdimen\figty

\newdimen\bwd
\bwd=2sp % 2sp (1/32768") is smallest visible width for Textures

\def\hline#1{\vbox{\smash{\hbox to #1{\leaders\hrule height \bwd\hfil}}}}

\def\vline#1{\hbox to 0pt{%
  \hss\vbox to #1{\leaders\vrule width \bwd\vfil}\hss}}

\def\clap#1{\hbox to 0pt{\hss#1\hss}}
\def\vclap#1{\vbox to 0pt{\offinterlineskip\vss#1\vss}}

\def\hstutter#1#2{\hbox{%
  \setbox0=\hbox{#1}%
  \hbox to #2\wd0{\leaders\box0\hfil}}}

\def\vstutter#1#2{\vbox{
  \setbox0=\vbox{\offinterlineskip #1}
  \dp0=0pt
  \vbox to #2\ht0{\leaders\box0\vfil}}}

\def\crosshairs#1#2{
  \dimen1=.002\drawx
  \dimen2=.002\drawy
  \ifdim\dimen1<\dimen2\dimen3\dimen1\else\dimen3\dimen2\fi
  \setbox1=\vclap{\vline{2\dimen3}}
  \setbox2=\clap{\hline{2\dimen3}}
  \setbox3=\hstutter{\kern\dimen1\box1}{4}
  \setbox4=\vstutter{\kern\dimen2\box2}{4}
  \setbox1=\vclap{\vline{4\dimen3}}
  \setbox2=\clap{\hline{4\dimen3}}
  \setbox5=\clap{\copy1\hstutter{\box3\kern\dimen1\box1}{6}}
  \setbox6=\vclap{\copy2\vstutter{\box4\kern\dimen2\box2}{6}}
  \setbox1=\vbox{\offinterlineskip\box5\box6}
  \smash{\vbox to #2{\hbox to #1{\hss\box1}\vss}}}

\def\boxgrid#1{\rlap{\vbox{\offinterlineskip
  \setbox0=\hline{\wd#1}
  \setbox1=\vline{\ht#1}
  \smash{\vbox to \ht#1{\offinterlineskip\copy0\vfil\box0}}
  \smash{\vbox{\hbox to \wd#1{\copy1\hfil\box1}}}}}}

\def\drawgrid#1#2{\vbox{\offinterlineskip
  \dimen0=\drawx
  \dimen1=\drawy
  \divide\dimen0 by 10
  \divide\dimen1 by 10
  \setbox0=\hline\drawx
  \setbox1=\vline\drawy
  \smash{\vbox{\offinterlineskip
    \copy0\vstutter{\kern\dimen1\box0}{10}}}
  \smash{\hbox{\copy1\hstutter{\kern\dimen0\box1}{10}}}}}

\def\figtext#1#2#3#4#5{
  \setbox\figtbox=\hbox{#5}
  \dp\figtbox=0pt
  \figtx=-#3\wd\figtbox \figty=-#4\ht\figtbox
  \advance\figtx by #1\drawx \advance\figty by #2\drawy
  \dimen0=\figtx \advance\dimen0 by\wd\figtbox \advance\dimen0 by-\drawx
  \ifdim\dimen0>\egap\global\egap=\dimen0\fi
  \dimen0=\figty \advance\dimen0 by\ht\figtbox \advance\dimen0 by-\drawy
  \ifdim\dimen0>\ngap\global\ngap=\dimen0\fi
  \dimen0=-\figtx
  \ifdim\dimen0>\wgap\global\wgap=\dimen0\fi
  \dimen0=-\figty
  \ifdim\dimen0>\sgap\global\sgap=\dimen0\fi
  \smash{\rlap{\vbox{\offinterlineskip
    \hbox{\hbox to \figtx{}\ifgridlines\boxgrid\figtbox\fi\box\figtbox}
    \vbox to \figty{}
    \ifgridlines\crosshairs{#1\drawx}{#2\drawy}\fi
    \kern 0pt}}}}

% macros to add space to text on specified sides

\def\hpad#1#2#3{\hbox{\kern #1\hbox{#3}\kern #2}}
\def\vpad#1#2#3{\setbox0=\hbox{#3}\dp0=0pt\vbox{\kern #1\box0\kern #2}}

% macro to give one text string the apparent height of another

% macro to center one text string over another

\def\stack#1#2#3{\vbox{\offinterlineskip
  \setbox2=\hbox{#2}
  \setbox3=\hbox{#3}
  \dimen0=\ifdim\wd2>\wd3\wd2\else\wd3\fi
  \hbox to \dimen0{\hss\box2\hss}
  \kern #1
  \hbox to \dimen0{\hss\box3\hss}}}

% macros to hide size of trailing exponents

\def\hexp#1{%
  \setbox0=\hbox{${}^{#1}$}%
  \hbox to .5\wd0{\box0\hss}}

%=========================================================================%
% macros for matrices and arrows
%=========================================================================%

% typical usage:
%   \rightarrowmat{2pt}{4pt}{d & bd \cr \!-c & 0 \cr 0 & -ac \cr}

\def\bmatrix#1#2{{\smallmath\left[\vcenter{\halign
  {&\kern#1\hfil$##\mathstrut$\kern#1\cr#2}}\right]}}

\def\rightarrowmat#1#2#3{
  \setbox1=\hbox{\kern#2$\bmatrix{#1}{#3}$\kern#2}
  \,\vbox{\offinterlineskip\hbox to\wd1{\hfil\copy1\hfil}
    \kern 3pt\hbox to\wd1{\rightarrowfill}}\,}

\def\leftarrowmat#1#2#3{
  \setbox1=\hbox{\kern#2$\bmatrix{#1}{#3}$\kern#2}
  \,\vbox{\offinterlineskip\hbox to\wd1{\hfil\copy1\hfil}
    \kern 3pt\hbox to\wd1{\leftarrowfill}}\,}

\def\rightarrowbox#1#2{
  \setbox1=\hbox{\kern#1\hbox{\smallmath #2}\kern#1}
  \,\vbox{\offinterlineskip\hbox to\wd1{\hfil\copy1\hfil}
    \kern 3pt\hbox to\wd1{\rightarrowfill}}\,}

\def\leftarrowbox#1#2{
  \setbox1=\hbox{\kern#1\hbox{\smallmath #2}\kern#1}
  \,\vbox{\offinterlineskip\hbox to\wd1{\hfil\copy1\hfil}
    \kern 3pt\hbox to\wd1{\leftarrowfill}}\,}

%=========================================================================%
% quire macros for preview mode and making booklets
%=========================================================================%

% \legalbooklet{20} makes a booklet from legal paper in landscape
% orientation, where "20" is the page count. To preview, give a negative
% pagecount. Either print using the legal duplex option on a modern laser
% printer, or struggle to simulate this effect manually. Bind using a long
% reach stapler.

% \preview squeezes two pages side by side in landscape orientation. It
% is not suitable for printing, but ideal for previewing on a two page
% monitor.

% \twoup squeezes two pages onto letter paper in landscape mode,
% suitable for printing.

% Each of these macros calls the file "quire.tex"

\def\bookletdims{
  \hsize=5.25truein
  \vsize=7truein
}

\def\legalbooklet#1{
  \input quire
  \bookletdims
  \htotal=7.0truein
  \vtotal=8.5truein
  % below computed from above
  \hoffset=\htotal
  \advance\hoffset by -\hsize
  \divide\hoffset by 2
  \voffset=\vtotal
  \advance\voffset by -\vsize
  \divide\voffset by 2
  \advance\voffset by -.0625truein
  \shhtotal=2\htotal
  % below doesn't need to change
  \horigin=0.0truein
  \vorigin=0.0truein
  \shstaplewidth=0.01pt
  \shstaplelength=0.66truein
  \shthickness=0pt
  \shoutline=0pt
  \shcrop=0pt
  \shvoffset=-1.0truein
  \ifnum#1>0\quire{#1}\else\qtwopages\fi
}

\def\preview{
  \input quire
  \bookletdims
  \hoffset=0.1truein
  \vtotal=8.5truein
  \shhtotal=14truein
  % below computed from above
  \voffset=\vtotal
  \advance\voffset by -\vsize
  \divide\voffset by 2
  \advance\voffset by -.0625truein
  \htotal=2\hoffset
  \advance\htotal by \hsize
  % below doesn't need to change
  \horigin=0.0truein
  \vorigin=0.0truein
  \shstaplewidth=0.5pt
  \shstaplelength=0.5\vtotal
  \shthickness=0pt
  \shoutline=0pt
  \shcrop=0pt
  \shvoffset=-1.0truein
  \qtwopages
}

\def\twoup{
  \input quire
  \hsize=4.79452truein % 5.25/1.095
  \vsize=7truein
  \vtotal=8.5truein
  \shhtotal=11truein
  % below computed from above
  \hoffset=-2\hsize
  \advance\hoffset by \shhtotal
  \divide\hoffset by 6
  \voffset=\vtotal
  \advance\voffset by -\vsize
  \divide\voffset by 2
  \advance\voffset by -12truept
  \htotal=2\hoffset
  \advance\htotal by \hsize
  % below doesn't need to change
  \horigin=0.0truein
  \vorigin=0.0truein
  \shstaplewidth=0.01pt
  \shstaplelength=0pt
  \shthickness=0pt
  \shoutline=0pt
  \shcrop=0pt
  \shvoffset=-1.0truein
  \qtwopages
}

%=========================================================================%
% timestamp (adapted from eplain.tex)
%=========================================================================%

\newcount\countA
\newcount\countB
\newcount\countC

\def\monthname{\begingroup
  \ifcase\number\month
    \or January\or February\or March\or April\or May\or June\or
    July\or August\or September\or October\or November\or December\fi
\endgroup}

\def\dayname{\begingroup
  \countA=\number\day
  \countB=\number\year
  \advance\countA by 0 % adjust after each leap day
  \advance\countA by \ifcase\month\or
    0\or 31\or 59\or 90\or 120\or 151\or
    181\or 212\or 243\or 273\or 304\or 334\fi
  \advance\countB by -1995
  \multiply\countB by 365
  \advance\countA by \countB
  \countB=\countA
  \divide\countB by 7
  \multiply\countB by 7
  \advance\countA by -\countB
  \advance\countA by 1
  \ifcase\countA\or Sunday\or Monday\or Tuesday\or Wednesday\or
    Thursday\or Friday\or Saturday\fi
\endgroup}

\def\timename{\begingroup
   \countA = \time
   \divide\countA by 60
   \countB = \countA
   \countC = \time
   \multiply\countA by 60
   \advance\countC by -\countA
   \ifnum\countC<10\toks1={0}\else\toks1={}\fi
   \ifnum\countB<12 \toks0={\sevenrm AM}
     \else\toks0={\sevenrm PM}\advance\countB by -12\fi
   \relax\ifnum\countB=0\countB=12\fi
   \hbox{\the\countB:\the\toks1 \the\countC \thinspace \the\toks0}
\endgroup}

\def\timestamp{\dayname, \the\day\ \monthname\ \the\year, \timename}

%==========================================================================
% macros (specific to this paper)
%==========================================================================

% surround with $ $ if not already in math mode
\def\enma#1{{\ifmmode#1\else$#1$\fi}}

% 
%\magnification \magstep1
%\showlabelsabove
\input diagrams.tex
%\vsize=5.0in
%\hsize=17truecm
\overfullrule=0pt
% Gothic fonts from AMSTeX
\font\tengoth=eufm10 Ê\font\fivegoth=eufm5
\font\sevengoth=eufm7
\newfam\gothfam Ê\scriptscriptfont\gothfam=\fivegoth
\textfont\gothfam=\tengoth \scriptfont\gothfam=\sevengoth

%
% Bold italic fonts
\font\tenbi=cmmib10 Ê\font\fivebi=cmmib5
\font\sevenbi=cmmib7
\newfam\bifam Ê\scriptscriptfont\bifam=\fivebi
\textfont\bifam=\tenbi \scriptfont\bifam=\sevenbi

\font\hd=cmbx10 scaled\magstep1
\def \fix#1 {{\hfill\break \bf (( #1 ))\hfill\break}}

\def\dim{\mathop{\rm dim}\nolimits}

\def\codim{\mathop{\rm codim}\nolimits}

\def\QED{{\hfill \qed\bigskip}}

\def\CC{{\bf C}}
\def\PP{{\bf P}}

\def\mm{{\bf m}}

\forward {degree equality}{Lemma}{1.2}
\forward {eghp}{Theorem}{7.1}

\centerline{\hd Row Ideals and Fibers of Morphisms}
\bigskip
\centerline {\bf David Eisenbud and Bernd Ulrich\footnote{\rm*}{\rm Both authors were supported in part by the NSF.
The second author is grateful to MSRI, where most of this research was done}
}\medskip
\centerline{\it Affectionately dedicated to Mel Hochster, who has been an inspiration to us for many years,}
\centerline{\it
on the occasion of his 65th birthday.}
\bigskip
\medskip

\noindent{\bf Abstract } We study the fibers of  projective
morphisms and rational maps.
We characterize the analytic spread of a homogeneous
ideal through properties of its syzygy matrix.
Powers of linearly presented ideals need not be linearly presented, but we identify a weaker linearity property that
is preserved by taking powers. 
\section{} Introduction

In this note we study the fibers of a rational map from an
algebraic point of view. We begin by describing four ideals related
to such a fiber.

Let $S=k[x_0,\dots,x_{n}]$ be a polynomial ring over an
infinite field $k$ with homogeneous maximal ideal $\mm$,
$I\subset S$ an ideal generated by
an $r+1$-dimensional vector space $W$ of forms
of the same degree, and $\phi$ the associated rational map $ \/ \PP^n\to \PP^r = \PP(W)$.
We will use this notation throughout. Since we are interested in the rational map,
we may remove common divisors of $W$, and thus
assume that $I$ has codimension at least 2.

A $k$-rational point $q$ in the target $\PP^r=\PP(W)$ is
by definition a codimension 1 subspace
$W_q$ of $W$. We write $I_q \subset S$ for the ideal generated by $W_q$. By a
homogeneous presentation of $I$ we will always mean a homogeneous free presentation
of $I$ with respect to a homogeneous minimal generating set.
If $F\to G=S\otimes W$ is such a presentation, then the composition
$F\to G \to S\otimes(W/W_q)$ is called the
{\it generalized row\/} corresponding to $q$, and its image is called the
{\it generalized row ideal\/} corresponding to $q$. It is the
ideal generated by the entries of a row in the homogeneous presentation matrix after
a change of basis. From this we
see that the generalized row ideal corresponding to $q$ is simply $I_q:I$.

The rational map $\phi$ is a morphism away from the algebraic
set $V(I)$, and we may form the fiber (=preimage) of the morphism over a point
$q\in \PP^r$. The saturated ideal of the scheme-theoretic
closure of this fiber is $I_q:I^\infty$, which we call the
{\it morphism fiber ideal\/} associated to $q$.

The rational map $\phi$ gives rise to a {\it correspondence\/}
$\Gamma\subset \PP^n\times \PP^r$, which is the closure of the
graph of the morphism induced by $\phi$. There are projections
$$
\PP^n\lTo^{\pi_1} \Gamma \rTo^{\pi_2} \PP^r
$$ and we define the {\it correspondence fiber\/} over $q$ to
be $\pi_1(\pi_2^{-1}(q))$. Since $\Gamma$ is
$\/ {\rm BiProj}( {\cal R})$, where ${\cal R}$ is the Rees algebra
$S[It]\subset S[t]$ of $I$, the correspondence fiber
is defined by the ideal $$
(I_qt {\cal R}:(It)^\infty)\cap S = \bigcup_i (I_qI^{i-1}:I^i).
$$
This ideal describes the locus where $I$ is not integral over $I_q$.
It is not hard to see that our four ideals are contained, each in the next,
$$\eqalign{
I_q \, &\subset \, I_q:I  \, \ \ \ \ \ \ \ \ \ \ \ \ \ \ \  \hbox{ row ideal } \cr
&\subset \, \bigcup_i (I_qI^{i-1}:I^i) \ \ \ \hbox{ correspondence fiber ideal } \cr
&\subset \, I_q:I^{\infty} \ \ \ \ \ \ \ \ \ \ \ \ \  \hbox{ morphism fiber ideal } .
}
$$
In Section 2 we compare the row ideals, morphism fiber ideals,
and correspondence fiber ideals.

In Section 3 we use generalized row ideals
 to give bounds on the analytic spread of $I$ by interpreting the analytic spread
as 1 plus the dimension of the image of $\phi$.

Many interesting rational maps $\phi$ are associated as above
to ideals $I$ with linear presentation matrices---see
for example Hulek, Katz and Schreyer [1992]. Thus
we are interested in linearly presented ideals and their
powers, which arise in the study of the graph. It is known
that the powers of a linearly presented ideal need not
be linearly presented. The first such examples were
exhibited by Sturmfels [2000]; for a survey of what is
known, see Eisenbud, Huneke and Ulrich [2006]. In Section 3
we also give criteria for birationality of the map, or for its restriction
to a linear subspace of $\PP^n$.

In Section 4 we generalize the notion of linear presentation
(of an ideal or module) in various directions:
A graded $S$-module $M$ generated by finitely many elements of the
same degree has {\it linear generalized row ideals\/}
if the entries of {\it every\/} generalized row of a homogeneous presentation
matrix for $M$
generate a linear ideal, i.e., an ideal generated by
linear forms. Obviously, any module with a linear presentation 
has this property, and we conjecture that the two notions are
equivalent in the case of ideals. The corresponding conjecture
is false for modules, but we prove it for modules of projective
dimension one. The main result of the section implies the
weak linearity property of powers mentioned in the abstract. It says, 
in particular, that if an ideal $I$ has linear generalized row ideals, 
then every power of $I$ has a homogeneous presentation
all of whose (ordinary) rows generate linear ideals.

\bigskip

\section{} Comparing the notions of fiber ideals

Recall that the row ideal for a point $q$ is always contained
in the Êcorrespondence fiber ideal, which is
contained in the morphism fiber ideal.
If the Êrow ideal is generated by linear forms (or, more
generally, is prime) and does
not contain $I$, then they are all equal. But in general
the containments are both strict:

\example{A} Let $S=k[a,b,c,d]$,
$J=(ab^2, ac^2, b^2c,bc^2)$, and $I=J+(bcd)$.
One can check that $I$ is linearly presented.
Computation shows that
the row ideal $J:I$ Êis $(b,c)$, while the correspondence
fiber ideal is $(a^2,b,c)$ and the morphism fiber ideal
is the unit ideal $J:I^\infty = S$. We
have no example of an $\mm$-primary ideal (regular morphism)
where all three are different: in the examples we have tried,
the correspondence fiber is equal to the morphism fiber.
(Of course for any regular map all three are equal up to saturation,
but we do not see why any two should be equal as ideals.)

\bigskip

Before stating the next result we recall that an ideal $I$ in a Noetherian ring
is said to be {\it of linear type\/} if the natural map from the symmetric
algebra of $I$ onto the Rees algebra of $I$ is an isomorphism.
If $I$ is of linear type, then $I$ cannot be
integral over any strictly smaller ideal, as can be seen by applying Theorem 4 
on p.152 of Northcott and Rees [1954] to the localizations of $I$.
We say that an ideal is {\it proper\/}
if it is not the unit ideal.

\proposition{linear corr ideals} If $I$ has linear generalized
row ideals, then every proper morphism fiber ideal is 
equal to the corresponding row ideal and hence generated by linear forms.
If $I$ is also of linear type on the punctured spectrum,
then  every proper correspondence fiber ideal is equal to the corresponding
row ideal.

\proof Suppose that the morphism fiber ideal
$I_{q}:I^\infty$
is not the unit ideal. In particular $I_{q}:I$ does not contain $I$.
The required equality for the first statement is
$$
I_{q}:I = I_{q}:I^\infty,
$$
which follows because $I_{q}:I$ is linear,
and thus prime.

Now suppose that $I$ is of linear type on the punctured spectrum, and that
the
correspondence fiber ideal
$H:=\bigcup_i Ê(I_{q}I^{i-1}:I^i)$ is proper. Set $K=I_q:I$, the row ideal.
We must show 
$K=H$. Since $K \subset H$ we may harmlessly assume that $K$ is not
$\mm$, the homogeneous maximal ideal of $S$.
By hypothesis the row ideal $K$ is generated by linear forms,
so it is prime. Since the localized ideals
$(I_{q})_K$ and $I_K$ are not equal, and $I_K$ is of
linear type, it follows that $I_K$ is not integral over
$(I_{q})_K$. Therefore
$H_K$ is a proper ideal. It follows that $H\subset K$,
as required. \QED

\example{} The last statement of
\ref{linear corr ideals} would be false without the
hypothesis that $I$ is of linear type on the punctured
spectrum. This is shown by \ref{A}.

\smallskip

\example{}
Let $Q$ be a quadratic form in $x_0,x_1,x_2$, and let
$F$ be a cubic form relatively prime to $Q$. The
rational map defined by $x_0Q,x_1Q,x_2Q, F$ has one
morphism fiber (and correspondence fiber) ideal $(Q)$, though for a general point
in the image both the morphism fiber ideal and the correspondence fiber ideal
are linear. This
example shows that in Theorem 4.1 of Simis [2004],
the point $p$ should be taken to be general.

\bigskip

\section{} How to compute the analytic spread
and test birationality

The notions of row ideals and fiber ideals provide
tests for the birationality of the map $\phi$ and lead to formulas for the analytic spread of the
ideal $I$. In our setting, the {\it analytic spread}
$\, \ell(I)$ of $I$ can be defined as one plus the dimension of the image
of the rational map $\phi$. Its ideal theoretic significance is
that it gives the smallest number of generators of a homogeneous ideal over which $I$ is integral, or equivalently, the smallest number of generators of an ideal in $S_{\mm}$ over which $I_{\mm}$
is integral, see the corollary on p.151 of Northcott and Rees [1954].

\medskip

\proposition{computing}
\item{$($a$)$} If $q$ is a point in $\, Ê\PP^r=\PP(W)$ such that $I_q : I^{\infty} \neq S$, then
$$
\ell(I) \geq 1+\codim (I_q :I^\infty) \, ;
$$
\item{$($b$)$} If $p$ is a general point in $\, \PP^n$, then
$$
\ell(I)=1+\codim (I_{{\phi}(p)}:I^\infty) \, ;
$$
\item{$($c$)$} If
there exits a point $q$ so that the row ideal $I_q:I$ is linear
of codimension $n$ and does not contain $I$, then 
$\phi$ is birational onto its image. Moreover, 
$\phi$ is birational onto its image if and only if
$I_{{\phi}(p)}:I^\infty$ is a linear ideal of codimension $n$ for a general point $p$.

\proof Set $J=I_{{\phi}(p)}$. If the ideal $I_q : I^{\infty}$ is proper it cannot be $\mm$-primary, and hence defines a non-empty fiber of the
morphism $\phi$. On the other hand, $J:I^\infty$ is
the defining ideal of a general fiber of the map. Thus
the dimension formula and the semicontinuity of fiber dimension, Corollary 14.5 and Theorem 14.8(a) in Eisenbud [1995], show that
$$
\codim (I_q :I^\infty) \leq \codim (J :I^\infty)= \dim \, Ê{\rm im} (\phi).
$$
However, the latter dimension is $\ell(I) -1$, proving parts
$(a)$ and $(b)$.

The second assertion in $(c)$ holds because the map 
is birational onto its image iff the general fiber is a reduced rational point.

We reduce the first assertion of $(c)$ to the second one. Assume that the row
ideal $I_q:I$ is linear of codimension $n$ and does not contain $I$.
Since $I_q:I$ is a prime ideal not containing $I$ it follows that
$I_q:I^{\infty} = I_q :I \neq S$. Thus the morphism fiber over $q$ is not
empty, and there exists a point $p \in \PP^n$ with $q=\phi(p)$.

Now let $T_0, \ldots, T_r$ be variables over $S$ and let $A_1$ denote 
the linear part of a homogeneous presentation matrix of $I$.
We can write
$(T_0,\dots,T_r)*A_1 = (x_0,\dots,x_n)*B$ for some matrix
$B$ whose entries are linear forms in the variables $T_i$
with constant coefficients. The dimension of the space
of linear forms in the row ideal corresponding to any
point $\phi(p)$ is the rank of $B$ when the coordinates of
$\phi(p)$ are substituted for the $T_i$; it is therefore
semicontinuous in $p$. Thus for $p$ general,
the dimension of the space of linear forms
in the ideal $I_{\phi(p)}:I$ is at least $n$, and then the same holds for
$J: I^{\infty}$. As this ideal defines a nonempty fiber, it
is indeed linear of codimension $n$.
\QED

Sometimes one can read off a lower bound on the analytic spread
even from a partial matrix of syzygies. The following result is inspired
by Proposition 1.2 of Hulek, Katz and Schreyer [1992].

\proposition{HKS}
With notation as above, suppose that $A$ is a  matrix of homogeneous
forms, each of whose columns is a syzygy on the generators of $I$. Let $A_q$ be the ideal generated by the elements of
the generalized row of $A$ corresponding to a point $q\in \PP^r$. If
there exists a prime ideal $P\in V(A_q)$ such that $A\otimes \kappa(P)$
has rank $r$,
then $I_q:I^\infty \neq S$ and
$$\ell(I)\geq 1+ \codim A_q \, .$$

\proof Since $A_q \subset I_q : I^{\infty}$, \ref{computing}$(a)$ Êshows
that the second claim follows from the first one.
To prove the first assertion,
$I_q:I^\infty \neq S$, it suffices to verify that
$(I_q:I^\infty)_P\neq S_P$.

As $A_P$ contains an $r\times r$ invertible submatrix, and these relations
express each generator of $I_P$ in terms of the one corresponding to $q$, it follows
that $A_P$ is a full presentation matrix of the ideal $I_P$. Thus
$(A_q)_P=(I_q:I)_P$. Furthermore, since $I_P$ is generated by
one element, and $I$ has codimension at least 2 by our blanket assumption,
it follows that $I_P=S_P$, whence $(A_q)_P=(I_q:I)_P=(I_q: I^{\infty})_P$.
On the other hand, $P\in V(A_q)$, so $(A_q)_P\neq S_P$, and we are done.
\QED

As in Proposition 1.2 of Hulek, Katz and Schreyer [1992], this gives criteria for birationality:

\corollary{} As in \ref{HKS} suppose that $A\otimes \kappa(P)$
has rank $r$ for some prime ideal $P\in V(A_q)$. The map $\phi$ is birational
onto its image if $A_q$ defines a reduced rational point in $\PP^n$.
The map $\phi$, restricted to a general $ \, Ê\PP^r\subset \PP^n$ is birational $($a Cremona
transformation$)$ if $A_q$ defines a reduced linear space of codimension $r$ in $\/ \PP^n$.

\proof Notice that $A_q \subset I_q:I \subset I_q : I^{\infty}$, where $I_q:I^{\infty} \neq S$
according to Proposition 3.2. Thus if $A_q$ defines a reduced rational point in $\PP^n$, then
the row ideal $I_q:I$ is linear of codimension $n$ and does not contain $I$. Thus $\phi$
is birational onto its image according to Proposition 3.1(c).

The second assertion follows from the first one, applied to the restriction of $\phi$.
\QED

For other, related criteria for birationality we refer to Simis [2004].

\bigskip

\section{} Ideals with linear row ideals and their powers

We begin this section by clarifying the relation between these properties of an ideal
or module:
to have a linear presentation matrix, to have linear generalized row ideals,
and to have {\it some} homogeneous presentation matrix all of whose row ideals are linear. Obviously,
if a presention matrix is linear then all its generalized row ideals are linear. However,
the converse does not hold, at least for the presentation of modules with torsion.
This can be seen by taking the matrix
$$
\pmatrix{
s&t&t^2\cr 0&s&0
}
$$
for instance. However, we have:

\proposition{} If $M$ is a graded $S$-module of projective dimension 1 
generated by finitely many homogeneous elements of the same degree,
and $M$ has linear generalized row ideals, then $M$ has a
linear presentation.

\proof Reduce modulo $n$ general linear forms, and use the Fundamental
Theorem for modules over principal ideal domains. \QED

Next, whenever an ideal has linear generalized row ideals, then obviously
there is a presentation matrix with only linear row ideals. Again, the
two concepts are not equivalent:

\example{} We consider the ideal $I=(s^4,s^3t,st^3,t^4) \subset S=\CC[s,t]$
corresponding to the morphism whose image is the smooth rational quartic curve
in $\PP^3$. A homogeneous presentation of this ideal is given by
$$
\diagram
S^2(-5)\oplus S(-6)&
\rTo^{
\pmatrix{
t& 0&0\cr
-s & 0 & t^2\cr
0&t&-s^2\cr
0&-s&0}
}
&
S^4(-4)
&
\rTo^{
\pmatrix{s^4&s^3t&st^3&t^4
}}
&
S
\enddiagram \ .
$$
The row ideals of the second and third rows in this
presentation are not linear. However, a change of
basis in $S^4(-4)$, corresponding to a different
choice of generators of $I$, makes them linear:
$$
\diagram
S^2(-5)\oplus S(-6)&
\rTo^{
\pmatrix{
t& 0&0\cr
0&s&0\cr
s-t & s-t & s^2-t^2\cr
-s+it&-is-t&s^2+t^2
}
}
&
S^4(-4)
&
\rTo^{
\pmatrix{F_0,\dots,F_3}
}
&
S
\enddiagram \ ,
$$
where
$$\eqalign{
F_0=&-s(s-t)(s^2+t^2+(s+t)(s-it))\cr
F_1=&-t(s-t)(s^2+t^2+(s+t)(is+t))\cr
F_2=&st(s^2+t^2)\cr
F_3=&-st(s^2-t^2) \ .
}
$$

\bigskip

Whereas powers of linearly presented ideals need not be
linearly presented,
the next result implies that having a homogeneous presentation with linear
row ideals is a weak linearity property that is indeed preserved when taking powers.

\theorem{linear} If $I$ has a homogeneous presentation matrix
where at least one row ideal is linear of codimension at least $\ell(I) -1$
and does not contain $I$, then each power of $I$ has some homogeneous presentation
matrix all of whose row ideals are linear of codimension
$\ell(I)-1$ and do not contain $I$.

\proof According to \ref{computing}$(b)$ for general $p \in \PP^n$, the morphism fiber ideal
$I_{{\phi}(p)}: I^{\infty}$ has codimension $\ell (I) -1$, and hence the row ideal $I_{{\phi}(p)}: I$
has codimension at most $\ell (I) -1$. Now one sees as in
the proof of \ref{computing}$(c)$ that Ê$I_{{\phi}(p)}: I$
is linear of codimension $\ell(I) -1$ and does not contain $I$.

Let $E=V(I)$ be the exceptional locus of $\phi$.
For each $d\geq 1$ the rational
map $\phi_d$ defined by the vector space of forms $W^d$
is regular on $\PP^n\setminus E$. 
For any point $p\in \PP^n \setminus E$, the ideal of $\phi(p) \in \PP(W)$
is generated by the vector space of linear forms $W_{\phi(p)}$, so the vector space of
forms of degree $d$ that it contains
is $W_{\phi(p)}W^{d-1}$. Thus
$(W^d)_{\phi_d(p)} = W_{\phi(p)}W^{d-1}$, and hence the row
ideal corresponding to $\phi_d(p)$ is 
$I_{\phi(p)}I^{d-1}:I^d$.

We now show that for general $p$, 
the row ideal $I_{\phi(p)}I^{d-1}:I^d$ is linear of codimension $\ell(I)-1$ and does 
not contain $I$. 
For trivial reasons we have
$$
I_{\phi(p)}:I\quad \subset \quad I_{\phi(p)}I^{d-1}:I^d
\quad \subset \quad I_{\phi(p)}I^{d-1}:I^\infty
\quad \subset \quad I_{\phi(p)}:I^\infty.
$$
By the above, 
$I_{\phi(p)}:I$ is a linear ideal of codimension $\ell (I) -1$ and does not contain $I$. Hence
$$
I_{\phi(p)}:I = I_{\phi(p)}:I^\infty,
$$
and therefore
$$
I_{\phi(p)}:I\quad = \quad I_{\phi(p)}I^{d-1}:I^d.
$$

Let $\dim W^d = N+1$.
Because the image of $\phi_d$ is nondegenerate,
$N+1$ general points of $\PP^n$
correspond to the $N+1$ rows of a
presentation matrix of $I^d$, so we are done.
\hfill \qed

\corollary{}
If $I$ has linear presentation, or even just linear generalized row ideals,
then every power of $I$ has a homogeneous presentation
matrix all of whose row ideals are linear of codimension
$\ell(I)-1$.

\proof According to \ref{computing}$(b)$, the homogeneous presentation matrix of
$I$ has a row ideal $I_q:I$
so that $\codim(I_q:I^{\infty}) = \ell(I) -1$. In particular $I_q:I^{\infty} \neq S$ and
hence $I$ is not contained in $I_q:I$. As $I_q:I$ is a linear ideal
we conclude that $I_q:I= I_q:I^{\infty}$, which gives $\codim(I_q:I) = \ell(I) -1$.
Now apply \ref{linear}.
\QED

\proposition{} Every ideal has a homogeneous presentation where every row ideal
has codimension at most $\ell(I)-1$

\proof Take a homogeneous presentation whose rows correspond to the
fibers through points of $
\, \PP^n$ not in the exceptional locus. The row ideals
are contained in
the morphism fiber ideals, which have codimension at most $\ell(I) - 1$ according
to Proposition 3.1(a). \QED

\bigskip

\section{} Some open problems

We would very much like to know the answer to the following questions:

\medskip

\item{1.} Can the homogeneous minimal presentation of an ideal $I$
have linear generalized row ideals without actually being
linear?

\item{2.} If $\phi$ is a regular map (that is, $I$ is $\mm$-primary),
are the correspondence
fiber ideals equal to the morphism fiber ideals? More generally,
when are the correspondence fiber ideals saturated with respect to $\mm$?

\item{3.} If $I$ is $\mm$-primary and linearly
presented, is every correspondence fiber ideal of
the morphism defined by $I^d$
either linear or $\mm$-primary?

\item{4.} Find lower bounds for
the number of linear relations $I^d$ could have
in terms of the number of linear relations on $I$. How close
can one come to the known examples?

\bigskip
\medskip

\references

D. Eisenbud: Commutative Algebra with a View Toward Algebraic Geometry.
Graduate Texts in Mathematics, 150. Springer-Verlag, New York, 1995.
\smallskip

D. Eisenbud, C. Huneke and B. Ulrich:
The regularity of Tor and graded Betti numbers.
Amer. J. Math. 128 (2006) 573--605.
\smallskip

K. Hulek, S. Katz and F.-O. Schreyer:
Cremona transformations and syzygies.
Math. Z. 209 (1992) 419--443.
\smallskip

D. G. Northcott and D. Rees: Reductions of ideals in local rings. Proc.
Cambridge Philos. Soc. 50 (1954) 145--158.
\smallskip

A. Simis:
Cremona transformations and some related algebras. J. Algebra 280 (2004) 162--179.

\smallskip
B. Sturmfels: Four counterexamples in combinatorial algebraic
geometry.
J. Algebra 230 (2000) 282--294.

\bigskip
\bigskip
\smallskip

\vbox{\noindent {\bf Author Addresses}\par
\smallskip
\noindent{David Eisenbud}\par
\noindent{Department of Mathematics, University of California, Berkeley,
Berkeley, CA 94720}\par
\noindent{eisenbud@math.berkeley.edu}\par
\medskip
\noindent{Bernd Ulrich}\par
\noindent{Department of Mathematics, Purdue University, West Lafayette, IN 47907}\par
\noindent{ulrich@math.purdue.edu}\par
}

\end